\documentclass[journal,10pt]{IEEEtran}
\usepackage{ifpdf}

\ifCLASSINFOpdf
\else
\fi
\usepackage{amsmath}
\usepackage{algorithmicx}
\usepackage{algpseudocode}

\usepackage{array}

\usepackage{amsfonts}
\usepackage{amssymb}
\usepackage{graphicx}
\usepackage{float}
\usepackage{setspace}

\begin{document}
%
\title{Partial estimators and application to covariance estimation of gaussian and elliptical distributions}


\author{\IEEEauthorblockN{Christophe Culan 
		,
		Claude Adnet\IEEEauthorrefmark{1}, 
	}
	\IEEEauthorblockA{ 
		Advanced Radar Concepts division,
		Thales, Limours 91338, France}
	\thanks{Manuscript received October XX, 2016; revised October XX, 2016. 
		Corresponding author: M. Culan (email: christophe.culan@thalesgroup.com).}}

\markboth{ArXiv prepublication}
{Culan \MakeLowercase{\textit{et al.}}: Partial estimators and application to covariance estimation of gaussian and elliptical distributions}
%



\IEEEtitleabstractindextext{%
\begin{abstract}
	Robustness to outliers is often a desirable property of statistical estimators. Indeed many well known estimators offer very good optimal performance in theory but are unusable in applied contexts because of their sensitivity to outliers. Of particular interest to the authors is the case of covariance estimators in adaptive matched filtering schemes in signal processing applications such as RADAR and SONAR detection, for which a contamination by outliers of the estimated noise covariance can lead to a great impact on performances, in particular when these outliers are similar to the target signal of the matched filter.\\
	This paper presents a generic method for building partial estimators from known estimators, which aim at avoiding these issues; the resulting algorithms are shown for a few chosen cases.
\end{abstract}

}

\maketitle

\IEEEdisplaynontitleabstractindextext

%
\IEEEpeerreviewmaketitle

\section{Introduction: the outlier issue in covariance estimations}
\label{sec:intro}

The issue of contamination by outliers is widely discussed in the statistical litterature \cite{daszykowski2007robust}\cite{pena2012multivariate}\cite{cousineau2015outliers}, and is of particular interest in adaptive filtering schemes in which they can greatly degrade performances.\\
Many have proposed selection schemes based on various definition of distances between signals in order to select samples which are most likely not to be outliers\cite{decurninge2014burg}\cite{decurninge2016robust}\cite{aubry2013median}. We propose in this paper to instead focus on the likelihood of the samples.

\subsection{Notations and conventions}
\label{ssec:notations}

In the following development, the following notations and conventions shall be observed:
\begin{itemize}
	\item For any topological space $\mathcal{R}_1$ and $\mathcal{R}_2$, the topological space $\mathcal{R}_1 \times \mathcal{R}_2$ is the product of $R_1$ and $R_2$, whereas $\mathcal{R}_1 \vee \mathcal{R}_2$ is the disjoint union of $R_1$ and $R_2$. One also notes, for a topological space $\mathcal{R}$,  $\mathcal{R}^N$ to be the product of $N$ copies of $\mathcal{R}$ and $\vee^N \mathcal{R}$ to be the disjoint union of $N$ copies of $\mathcal{R}$.
	\item Let $\mathcal{R}_1$, ..., $\mathcal{R}_N$ be topological spaces, and let $\mu_n$ be a measure defined on $\mathcal{R}_n$ for $1 \leq n \leq N$. The measure $\mu_1 \dots \mu_n = \prod_{n=1}^N \mu_n$ is the product measure of $(\mu_n)_{1\leq n\leq N}$ defined on $\prod_{n=1}^N \mathcal{R}_n$, whereas\\ $\mu_1 \vee \dots \vee \mu_n = \bigvee_{n=1}^N \mu_n$ is the measure of $\bigvee_{n=1}^N \mathcal{R}_n$ such that its restriction to $\mathcal{R}_n$ is $\mu_n$ for $1 \leq n \leq N$ and is the joint measure of $(\mu_n)_{1 \leq n \leq N}$.\\
	Moreover for a topological space $\mathcal{R}$ and a measure $\mu$ defined on $\mathcal{R}$, $\mu^N$ is the product measure of $N$ copies of $\mu$ defined on $\mathcal{R}^N$, whereas $\vee^N \mu$ is the joint measure of $N$ copies of $\mu$ defined on $N\cdot\mathcal{R}$.
	\item The vector space $\mathbb{C}^d$ is canonically identified to $\mathbb{R}^{2d}$ \cite{scharf1991statistical}.
	\item $\mathcal{H}_d(\mathbb{R})$ is the set of symmetric matrices of size $(d,d)$; $\mathcal{H}_d(\mathbb{C})$ is the set of hermitian matrices of size $(d,d)$.\\
	$\mathcal{H}_d^+(\mathbb{K})$ is the subset of matrices of 	$\mathcal{H}_d^+(\mathbb{K})$ which are positive; $\mathcal{HP}_d^+(\mathbb{K})$ is the subset of matrices of $\mathcal{H}_d^+(\mathbb{K})$ of unit determinant.
	\item The vector space $\mathbb{C}^d$ is canonically identified to $\mathbb{R}^{2d}$ \cite{scharf1991statistical}.
	\item  $X^\dagger$ is the transpose conjugate of any matrix or vector $X$.
	\item $\mu_{\mathbb{R}^d}$ is the Lebesgues measure of the space $\mathbb{R}^d$. Similarly, $\mu_{\mathbb{C}^d}$ is the Lebesgues measure of $\mathbb{C}^d \simeq \mathbb{R}^{2d}$.
	\item $\mathcal{S}_{d-1}$ is the $(d-1)$-sphere, which is identified as the following part of $\mathbb{R}^d$: $\mathcal{S}_{d-1} \simeq \left\{ x \in \mathbb{R}^d; x^\dagger x = 1 \right\}$. Similarly $\mathcal{S}_{2d-1}$ is identified to $\left\{ x \in \mathbb{C}^d; x^\dagger x = 1 \right\}$ in $\mathbb{C}^d$.
	\item $s_{d-1}$ shall denote the probability distribution of $\mathcal{S}_{d-1}$ isotropic for the canonical scalar product of $\mathbb{C}^d$, defined by:
	\begin{equation}
	\left<y,x\right> = \sum\limits_{k=0}^{d-1} \overline{y_k} x_k
	\end{equation}
	\item $\delta_x^\mathcal{R}$ is the Dirac delta distribution centered in $x$ in space $\mathcal{R}$.
	\item One shall note:
	\begin{equation}
	h_\delta = H\left(\delta_0^{\mathbb{R}}|\mu_{\mathbb{R}}\right) = \int_{t \in \mathbb{R}} \log\left(\frac{d\delta_0^{\mathbb{R}}}{dt}(t)\right)d\delta_0^{\mathbb{R}}(t)
	\end{equation}
	 It is the entropy of the Dirac distribution relative to the Lebesgues measure of $\mathbb{R}$ due to its translational symmetry, which is positive and infinite. It can be used to express different other entropies in higher dimensional settings:
	\begin{equation}
	\begin{array}{l}
	H\left(\delta_x^{\mathbb{R}^d}|\mu_{\mathbb{R}^d}\right) = d h_\delta\vspace{1.5 mm}\\
	H\left(\delta_x^{\mathbb{C}^d}|\mu_{\mathbb{C}^d}\right) = 2d h_\delta
	\end{array}
	\end{equation}
	\item We shall call a reference measure $\mu$ of a space $\mathcal{R}$ a measure such that $h_\delta(x) = H\left(\delta_x^\mathcal{R}|\mu\right)$ does not depend on $x$. Such a measure need not exist as there might not be any way of deciding of the relative order of two different values $h_\delta(x)$ and $h_\delta(y)$. This in general depends on the symmetries of the considered space. 
	One can note, for example, that $\mu_{\mathbb{R}^d}$ is a reference measure of $\mathbb{R}^d$, whereas $\mathbb{C}P^{d-1}$ and $\mathbb{R}P^{d-1}$ have no reference measures.
\end{itemize}

\section{General background of partial estimators}
\label{sec:general_partial}

Suppose that we have a family of distributions $\left(P_\theta\right)_{\theta \in \Theta}$, and a set of $N$ samples which consist of at least $pN$ i.i.d samples drawn from a distribution $P_{\theta_0}$; the remaining samples are of unknown characteristics and might follow the same distribution, or might be outliers. One then wants to build an estimator which can efficiently reject these outliers.

\subsection{Partial likelihood score}
\label{ssec:partial_likelihood}

As is shown in \cite{CVC_BT}, the likelihood score of a model distribution $P$ for a sampling distribution $S$ can be expressed as an absolute quantity, in terms of a relative entropy between the sampling distribution and the model distribution \cite{kullback1951information}\cite{akaike1998information}:

\begin{equation}
\label{eq:def_likelihood}
l(P|S) = -H(S|P) = \int_{x \in \mathcal{R}} \log\left(\frac{dP}{dS}(x)\right)dS(x)
\end{equation}  

In the case in which one knowns that some of the samples are outliers, one could derive a model based on a mixture of a true distribution model and an outlier model. However by definition,there is often no information on the distribution of the outliers. Instead we propose to define a likelihood conditional only on the value of the samples which are supposed to follow the true distribution.\\

Hence we give the following definition of the partial likelihood measure of a partial domain $X$ by:

\begin{equation}
\label{eq:def_partial_ondomain}
l_{\left| X\right.}(P|S) = -H_{\left| X\right.}(P|S) = \int_{x \in X} \log\left(\frac{dP}{dS}(x)\right)dS(x)
\end{equation}

The principle of maximum partial likelihood of order $p \leq 1$ over a distribution family $\left(P_\theta\right)_{\theta \in \Theta}$ can then be stated as:\\

Find $\theta \in \Theta,X \subset \mathcal{R}$ maximizing $l_{\left| X\right.}(P_\theta|S)$, subject to the constrain $S(X) \geq p$.\\

This is a direct generalization of the maximum likelihood estimation procedure. This in general can be solved by considering the concentrated version already maximized over partial domains $X$ such that $S(X) \geq p$:

\begin{equation}
\label{eq:def_partial_likelihood}
l_p(P|S) = -H_p(P|S) = \sup_{S(X) \geq p} \int_{x \in X} \log\left(\frac{dP}{dS}(x)\right)dS(x)
\end{equation}

\subsection{Maximization over the partial support for i.i.d samples}

In the following development, the sampling distribution is taken to be the standard joint of Dirac distributions over a base space $\mathcal{B}$; the model distribution is the same on each component \cite{CVC_BT}:  

\begin{equation}
\left\{\begin{array}{l}
S = \frac{1}{N}\bigvee_{n=1}^N \delta_{x_n}^{\mathcal{B}}\\
\mathcal{P}_\mathcal{R} = \vee^N \mathcal{P}
\end{array}\right.
\end{equation}

Moreover the space $\mathcal{B}$ is supposed to be of finite dimension $d$ over $\mathbb{R}$ and admits a reference volume $\mu$. The likelihood score of a distribution $P$ for the sampling distribution $S$ can then be re-expressed as \cite{CVC_BT}:

\begin{equation}
l(P|S) =  - dh_\delta + \sum\limits_{n=1}^N  \log\left(\frac{dP}{d\mu}(x_n)\right)
\end{equation}

with $h_\delta = H\left(\delta_0^\mathbb{R}|\mu_\mathbb{R}\right)$ being the entropy of the standard one-dimensional Dirac distribution relative to the reference measure $\mu$ of $\mathcal{B}$.\\	

The corresponding partial likelihood measure is given by:

\begin{equation}
l(P|S)= -\frac{\sharp\left(\left\{n; x_n \in X\right\}\right)}{N}d  h_\delta+\sum\limits_{k \in \left\{n;x_n \in X\right\}} \log\left(\frac{dP}{d\mu}(x_k)\right)
\end{equation}

Since $h_\delta$ is infinite, a maximizer of this partial likelihood necessarily minimizes $\sharp\left(\left\{n; x_n \in X\right\}\right)$, which is then equal to $\lceil pN \rceil$. Thus the maximization of the partial likelihood for order $p$ for a given model distribution $P$ is equivalent to finding a subfamily of indices $\mathcal{K} \in \mathcal{P}_{\lceil pN \rceil}\left(\left[\left[1;N\right]\right]\right)$ maximizing the corresponding likelihood score:

\begin{equation}
l_\mathcal{K}(P|S) = -\frac{\lceil pN \rceil}{N}(dh_\delta) + \sum\limits_{k \in \mathcal{K}} \log\left(\frac{dP}{d\mu}(x_k)\right)
\end{equation}

This can be done by finding a permutation $o$ ordering the values of the density with respect to the reference measure in increasing order:
\begin{equation}
\frac{dP}{d\mu}\left(x_{o(1)}\right) \geq \frac{dP}{d\mu}\left(x_{o(2)}\right) \geq ... \geq \frac{dP}{d\mu}\left(x_{o(N)}\right)
\end{equation}

The corresponding maximum value is then given by:

\begin{equation}
l_p(P|S) = -\frac{\lceil pN \rceil}{N}(dh_\delta) + \sum\limits_{n=1}^{\lceil pN \rceil} \log\left(\frac{dP}{d\mu}\left(x_{o(n)}\right)\right)
\end{equation}

Hence this corresponds to a selection of the $\lceil pN \rceil$ most likely samples in the estimation procedure.

\subsection{Partial estimators}

A complete optimization over a distribution family $\left(P_\theta\right)_{\theta \in \Theta}$ would then require the maximization of the concentrated likelihood $l_p(P|S)$. However this is a difficult problem due to the mixing of two types of optimization problems:

\begin{itemize}
	\item the optimization over the partial domain is a discrete optimization problem.
	\item the optimization over the parameter space $\Theta$ generally is a continuous optimization problem, which is quite often solved by numerical methods.
\end{itemize}

One could note however that if there is a way to compute a maximum likelihood estimator $e(S)$ associated to the distribution family $\left(P_\theta\right)_{\theta \in \Theta}$ for any given standard sampling distribution $S = \frac{1}{N}\bigvee_{n=1}^N \delta_{x_n}^\mathcal{R}$, one can in theory find a true maximum partial likelihood parameter of any order $p$: 

\begin{equation}
	e_p(S) = e(S_{\mathcal{K}_\text{max}})
\end{equation}

with:

\begin{equation}
\mathcal{K}_\text{max} = \text{argmax}_{\mathcal{K} \in \mathcal{P}_{\lceil pN \rceil}\left(\left[\left[1;N\right]\right]\right)} l(e(S_\mathcal{K})|S_\mathcal{K})
\end{equation}

with $S_\mathcal{K} = \frac{1}{\lceil pN \rceil} \bigvee_{k \in \mathcal{K}} \delta_{x_k}^\mathcal{R}$

However this is in general impractical; indeed one has:
\begin{displaymath} 
\sharp\left(\mathcal{P}_{\lceil pN \rceil}\left(\left[\left[1;N\right]\right]\right)\right) = \left(\begin{array}{c}
N\\
\lceil pN \rceil
\end{array}\right)
\end{displaymath}
Thus this quickly makes the computation intractable as $N$ grows, as one needs to compute $e(\mathcal{K})$ for each $\mathcal{K} \in \mathcal{P}_{\lceil pN \rceil}\left(\left[\left[1;N\right]\right]\right)$.\\

Instead we propose to resort to the following suboptimal partial estimation procedure:
supposing we already have an estimator $e(S)$ for any standard sampling distribution $S$ (not necessarily a maximum likelihood estimator), the following procedure can be used:

\begin{onehalfspace}
\begin{samepage}
\begin{algorithmic}[1]
	\Function{partial}{$e,\left(P_\theta\right)_{\theta \in \Theta},\theta_0,\left(x_n\right)_{1 \leq n \leq N},p,K_\text{max}$}
		\State $\theta \gets \theta_0$
		\For{$n$ \textbf{from} $1$ \textbf{to} $N$}
			\State $l_n \gets l\left(P_\theta|\delta_{x_n}^\mathcal{R}\right)$
		\EndFor\vspace{1.5 mm}
		\State $o_0 \gets \text{argsort}_{\downarrow}\left(\left(l_n\right)_{1 \leq n \leq N}\right)$\vspace{1.5 mm}
		\For{$k$ \textbf{from} $1$ \textbf{to} $K_\text{max}$}
			\State $\theta \gets e\left(S_{\left\{o_0(n);1 \leq n \leq \lceil pN \rceil\right\}}\right)$
			\For{$n$ \textbf{from} $1$ \textbf{to} $N$}
				\State $l_n \gets l\left(P_\theta|\delta_{x_n}^\mathcal{R}\right)$
			\EndFor\vspace{1.5 mm}
			\State $o_1 \gets \text{argsort}_{\downarrow}\left(\left(l_n\right)_{1 \leq n \leq N}\right)$\vspace{1.5 mm}
			\If{$o_0 = o_1$}
				\State \textbf{break}
			\Else

				\State $o_0 \gets o_1$
			\EndIf
		\EndFor
		\State \Return $\theta$
	\EndFunction\vspace{2 mm}
\end{algorithmic}
\end{samepage}
\end{onehalfspace}

Unfortunately this optimization procedure offers no guaranty of convergence, thus this should be checked empirically every time one desires to use it. \\
Moreover even when convergence is insured and a maximum likelihood estimator is used as the basis estimate, one should keep in mind that there is no guaranty of true maximization of the partial likelihood; thus the quality of the estimate should also be checked empirically.\\
Of final note is the fact that applying a partial estimation procedure on an otherwise unbiased estimator can produce a bias; such biases should be characterized whenever possible.

\section{Application to some covariance estimators}

This section discusses the application of the partial estimation procedure to different covariance estimators. The first part treats the case of multivariate gaussian variables, whereas the second part treats the case of robust Tyler-type estimators for elliptical distributions, ans is based on the development given in \cite{CVC_BT}. The third part treats the case of other M-estimators, and is based on the development of such estimators as likelihood maximizers given in \cite{CVC_BT}.

\subsection{Estimators under gaussian background hypothesis}

The likelihood of a delta distribution $\delta_x^{\mathbb{K}^d}$ or for a centered multivariate normal model of covariance $\Sigma$ is given, up to constant additive terms, by \cite{anderson1985maximum}:

\begin{equation}
l(\Sigma|\delta_{[x]}) = -\frac{c_\mathbb{K}}{2}\left( \log|\Sigma| + x^\dagger \Sigma^{-1} x\right)
\end{equation}

with $c_\mathbb{K}$ for $\mathbb{K} = \mathbb{R}$ and $c_\mathbb{K}=2$ for $\mathbb{K} = \mathbb{C}$. Note that no phase symmetry of the sampling distribution is taken into account in the complex case here \cite{CVC_BT}.

This likelihood is a strictly decreasing function of $x^\dagger \Sigma^{-1} x$; hence the ordering of the likelihoods and selection of the most likely samples can be done by selecting the samples for which the value of $x^\dagger \Sigma^{-1} x$ is the smallest. Note moreover that $\Sigma^{-1}$ need only be known up to a multiplicative constant.\\

Let us now see how this can be applied to some algorithms.

\subsubsection{Partial sample covariance}

The application of the partial estimation procedure to the sample covariance estimator can be used to obtain a robust estimator in this case. The corresponding procedure is outlined below:\\

\begin{onehalfspace}
\begin{samepage}
\begin{algorithmic}[1]
	\Function{partial\_SCM}{$\left(x_n\right)_{1 \leq n \leq N},K_\text{max}$}
		\State $\Sigma \gets I$
		\For{$n$ \textbf{from} $1$ \textbf{to} $N$}
			\State $\tau_n \gets {x_n}^\dagger \Sigma^{-1} x_n$
		\EndFor\vspace{1.5 mm}
		\State $o_0 \gets \text{argsort}_{\uparrow}\left(\left(x_n\right)_{1 \leq n \leq N}\right)$\vspace{1.5 mm}
		\For{$k$ \textbf{from} $1$ \textbf{to} $K_\text{max}$}\vspace{1.5 mm}
			\State $\Sigma \gets \frac{1}{\lceil pN \rceil}\sum\limits_{n=1}^{\lceil pN \rceil} x_{o_0(n)} {x_{o_0(n)}}^\dagger$\vspace{1.5 mm}
			\For{$n$ \textbf{from} $1$ \textbf{to} $N$}
				\State $\tau_n \gets {x_n}^\dagger \Sigma^{-1} x_n$
			\EndFor\vspace{1.5 mm}
			\State $o_1 \gets \text{argsort}_{\uparrow}\left(\left(\tau_n\right)_{1 \leq n \leq N}\right)$\vspace{1.5 mm}
			\If{$o_0 = o_1$}
				\State \textbf{break}
			\Else
				\State $o_0 \gets o_1$
			\EndIf
		\EndFor
		\State \Return $\Sigma$
	\EndFunction
\end{algorithmic}
\end{samepage}
\end{onehalfspace}
\vspace{2 mm}

One should note that this partial algorithm introduces a scale bias on the estimated covariance matrix; indeed the samples of largest scale are rejected from the estimation. Fortunately assuming convergence of the algorithm towards a true partial likelihood maximum, this bias can be expressed under the hypothesis that no outlier is present and depends only on the dimension of the space:

\begin{equation}
b(d) = \gamma\left(c_\mathbb{K}\frac{d}{2}+1,\gamma^{-1}\left(c_\mathbb{K}\frac{d}{2},\frac{\lceil pN \rceil}{N}\right)\right)
\end{equation}

with $\gamma$ being the incomplete gamma function and $\gamma^{-1}$ the incomplete inverse gamma function.

\subsubsection{Estimation under toeplitz constrain in the complex case: the partial multisegment Burg algorithm}

As noted in \cite{CVC_BT}, the Toeplitz constrain is already very difficult to enforce on a maximum likelihood estimation procedure. Instead we propose to base our estimate on the multisegment Burg estimator, which shall be noted as Burg($\left(x_n\right)_{1 \leq n \leq N}$) \cite{haykin1982maximum}\cite{burg1978maximum}\cite{ulrych1976time}. 
The corresponding partial estimator is given by:\\

\begin{onehalfspace}
\begin{samepage}
\begin{algorithmic}[1]
	\Function{Partial\_Burg}{$\left(x_n\right)_{1 \leq n \leq N}, p, K_\text{max}$}
		\State $\sigma^2 \gets 1$
		\State $\mu \gets (0)_{1 \leq m \leq d-1}$\vspace{1.5 mm}
		\State $\Sigma_{-1} \gets \textproc{Trench}(\sigma^2,\mu)$ 
		\For{$n$ \textbf{from} $1$ \textbf{to} $N$}
			\State $\tau_n \gets {x_n}^\dagger \Sigma_{-1} x_n$
		\EndFor\vspace{1.5 mm}
		\State $o_0 \gets \text{argsort}_{\uparrow}\left(\left(x_n\right)_{1 \leq n \leq N}\right)$\vspace{1.5 mm}
		\For{$k$ \textbf{from} $1$ \textbf{to} $K_\text{max}$}\vspace{1.5 mm}
			\State $\left(\sigma^2,\mu\right) \gets \textproc{Burg}\left(\left(x_{o_0(n)}\right)_{1 \leq n \leq \lceil pN \rceil}\right)$\vspace{1.5 mm}
			\State $\Sigma_{-1} \gets \textproc{Trench}(\sigma^2,\mu)$
			\For{$n$ \textbf{from} $1$ \textbf{to} $N$}
				\State $\tau_n \gets {x_n}^\dagger \Sigma_{-1} x_n$
			\EndFor\vspace{1.5 mm}
			\State $\sigma^2 \gets \frac{1}{d \lceil pN \rceil} \sum\limits_{n=1}^{\lceil pN \rceil} \tau_{o_0(n)}$\vspace{1.5 mm}
			\State $o_1 \gets \text{argsort}_{\uparrow}\left(\left(\tau_n\right)_{1 \leq n \leq N}\right)$\vspace{1.5 mm}
			\If{$o_0 = o_1$}
				\State \textbf{break}
			\Else
			\textbf{else do:}
				\State $o_0 \gets o_1$
			\EndIf
		\EndFor
		\State \Return $\Sigma_{-1}$
	\EndFunction
\end{algorithmic}
\end{samepage}
\end{onehalfspace}
\vspace{2 mm}

The \textproc{Burg} and \textproc{Trench} functions correspond respectively to the multisegment Burg algorithm  returning the residual error $\sigma^2$ and the Schur coefficients $(\mu_m)_{1 \leq m \leq d-1}$ \cite{haykin1982maximum}, and the Trench algorithm returning the inverse covariance matrix \cite{trench1964algorithm}\cite{zohar1969toeplitz}\cite{barbaresco1997analyse}.

Again this partial version of the algorithm suffers from a scale bias; under the hypothesis that no outliers are present in the sampling distribution, one can use $b(2d) =  \gamma\left(d+1,\gamma^{-1}\left(d,\frac{\lceil pN \rceil}{N}\right)\right)$ to approximate this bias.

\subsection{Partial estimators for elliptical distributions}

We now shift our study to the more general class of elliptical models \cite{ollila2012complex}\cite{chmielewski1981elliptically}. These elliptical models offer a broad generalization of the multivariate gaussian model which are used in statistical finance for portfolio modeling, as well as for the modeling of arbitrary impulsive distributions in signal processing applications, such as the clutter distribution in radar detection applications \cite{conte1995asymptotically}. Many widely used distribution models belong to this class such as K-distributions \cite{conte1991modelling}, t-distributions \cite{krishnaiah1986complex} \cite{ollila2003robust}, or the class of compound gaussian models which are widely used in simulations because they are easy to generate \cite{conte1987characterisation} \cite{gini2002vector} \cite{pascal2008covariance}.\\

Such distributions are of the form \cite{CVC_BT}:
\begin{equation}
dP(x) = dQ\left(\sqrt{x^\dagger R^{-1}x}\right)ds_{c_\mathbb{K}d-1}\left( \frac{L(R)^{-1} x}{\sqrt{x^\dagger R^{-1} x}} \right)
\end{equation}

with the correlation matrix $R$ being positive definite of unit determinant, and the radial distribution $Q$ verifying that $Q\left(\mathbb{R}_+\right)=1$.

In a manner similar to the treatment given in \cite{CVC_BT} to obtain the concentrated likelihood on the radial distribution, one can equivalently obtain the concentrated partial likelihood of a given order over the radial distribution and the partial domain for a standard sampling distribution $S = \frac{1}{N} \bigvee_{n=1}^N \delta_{[x_n]}^{[\mathbb{K}^d]}$:

\begin{equation}
\begin{split}
l_p(R|S) = &-\frac{\lceil pN \rceil}{N}\left(c_\mathbb{K} (d-1)h_\delta + \log\left(\frac{\Gamma\left(c_\mathbb{K}\frac{d}{2}\right)}{2\pi^{1+c_\mathbb{K}\left(\frac{d}{2}-1\right)}}\right)\right)\\
&-c_\mathbb{K} \frac{d-1}{2\lceil pN \rceil}\sum\limits_{n=1}^{\lceil pN \rceil} \log\left({x_n}^\dagger R^{-1} x_n\right)
\end{split}
\end{equation}

with $c_\mathbb{K}= 1$ for $\mathbb{K} = \mathbb{R}$ and $c_\mathbb{K}=2$ for $\mathbb{K} = \mathbb{C}$.\\

It follows that the expression of the concentrated likelihood for a single sample can be used for the data selection in order to maximize the likelihood over partial domains; this likelihood is given, up to constant terms, by:

\begin{equation}
l\left(R|\delta_{[x]}^{[\mathbb{K}^d]}\right) = - c_\mathbb{K} \frac{d-1}{2}\log\left(x^\dagger R^{-1} x\right)
\end{equation}

Interestingly this is again a strictly decreasing function of $x^\dagger R^{-1} x$; therefore the selection of the most likely samples can be done by keeping the $\lceil pN \rceil$ samples for which $x^\dagger R^{-1} x$ is smallest.

\subsubsection{partial Tyler fixed point algorithm}

Since Tyler's estimator is already obtained by a fixed point algorithm, we propose the following algorithm which incorporates the partial estimation directly in its fixed point loop \cite{tyler1987statistical}\cite{pascal2008covariance}:\\

\begin{onehalfspace}
\begin{samepage}
\begin{algorithmic}[1]
	\Function{pTyler}{$\left(x_n\right)_{1 \leq n \leq N}, p, \epsilon, K_\text{max}$}
		\State $R_{-1} \gets I$
		\For{$k$ \textbf{from} $1$ \textbf{to} $K_\text{max}$}
			\For{$n$ \textbf{from} $1$ \textbf{to} $N$}
				\State $\tau_n \gets {x_n}^\dagger R_{-1} x_n$
			\EndFor\vspace{1.5 mm}
			\State $o \gets \text{argsort}_{\uparrow}\left(\left(\tau_n\right)_{1 \leq n \leq N}\right)$\vspace{1.5 mm}
			\State $S \gets \sum\limits_{n=1}^{\lceil pN \rceil} \frac{x_{o(n)} {x_{o(n)}}^\dagger}{\tau_{o_{(n)}}}$\vspace{1.5 mm}
			\State $R \gets \frac{S}{\text{tr}(S)}$\vspace{1.5 mm}
			\If{$\text{tr}\left(\left(R_{-1} R - I\right)^2\right) \leq \epsilon$}\vspace{1.5 mm}
				\State \textbf{break}	
			\Else	
				\State $R_{-1} \gets R^{-1}$
			\EndIf
		\EndFor
		\State \Return $R_{-1}$
	\EndFunction
\end{algorithmic}
\end{samepage}
\end{onehalfspace}
\vspace{2 mm}

\subsubsection{Estimation under complex Toeplitz constrain : partial Burg-Tyler algorithm}

Similarly to Tyler's estimator, the \textproc{pTyler} algorithm can be adapted to take into account constrains on the covariance structure \cite{CVC_BT}. We shall consider the example of a Toeplitz constrain.\\
Since the problem of maximizing the likelihood under Toeplitz constrain is again quite hard, we propose to use the partial Burg-tyler algorithm, based on the Burg-Tyler algorithm introduced in \cite{CVC_BT}:

\begin{onehalfspace}
\begin{samepage}
\begin{algorithmic}[1]
	\Function{pBT}{$\left(x_n\right)_{1 \leq n \leq N}, p, \epsilon, K_\text{max}$}
		\State $\mu \gets (0)_{1 \leq m \leq d-1}$
		\For{$k$ \textbf{from} 1 \textbf{to} $K_{\max}$}
			\State $R_{-1} \gets \textproc{Trench}(1,\mu)$ 
			\For{$n$ \textbf{from} 1 \textbf{to} $N$}
				\State $\tau_n \gets {x_n}^\dagger R_{-1} x_n$
			\EndFor\vspace{1.5 mm}
			\State $o \gets \text{argsort}_{\uparrow}\left(\left(\tau_n\right)_{1 \leq n \leq N}\right)$\vspace{1.5 mm}
			\State $\left(\sigma^2,\nu\right) \gets \textproc{Burg}\left(\left(\frac{x_{o(n)}}{\sqrt{\tau_{o(n)}}}\right)_{1 \leq n \leq \lceil pN \rceil}\right)$\vspace{1.5 mm}
			\If{$\sum\limits_{m=1}^{d-1} (d-m) \text{atanh}^2\left(\left|\frac{\nu_m - \mu_m}{1-\overline{\mu_m} \nu_m}\right|\right) \leq \epsilon$}\vspace{1.5 mm} 
				\State \textbf{break}
			\Else
				\State $\mu \gets \nu$
			\EndIf
		\EndFor
		\State \Return $R_{-1}$
	\EndFunction
\end{algorithmic}
\end{samepage}
\end{onehalfspace}
\vspace{2 mm}

This algorithm is not per say a partial likelihood maximizer; however it is sufficient to outperform the \textproc{pTyler} algorithm on scenarii involving stationary signals, in all the tests performed by the authors so far on both simulated and real data.

\subsection{Extension to other M-estimators}

Such estimation procedures can also be extended to other M-estimators \cite{marona1998robust}\cite{ollila2003robust}. Indeed restricting oneself to elliptical distributions with radial distributions of the form $dQ(x) = \sigma dQ_0(\sigma x)$ with $\sigma >0$ and $Q_0$ being a fixed probability distribution, one can express the associated likelihood for a single sample in the following form \cite{CVC_BT}: 

\begin{equation}
l(Q,R|\frac{1}{N}\bigvee_{n=1}^N\delta_{[x_0]}^{\mathbb{K}^d}) = 
-\frac{c_\mathbb{K}}{2}\log\left|\Sigma\right|+\frac{1}{2} g\left({x_0}^\dagger \Sigma^{-1}x_0\right)
\end{equation}

with $\Sigma = \sigma R$. Thus this leads to the following partial likelihood:

\begin{flalign*}
l_p\left(Q,R|\frac{1}{N}\bigvee_{n=1}^N\delta_{[x_n]}^{\mathbb{K}^d}\right) = &&
\end{flalign*}
\begin{equation} 
-\frac{\lceil pN\rceil}{N}\frac{c_\mathbb{K}}{2}\log\left|\Sigma\right|+\frac{1}{2\lceil pN \rceil} \sum\limits_{n=1}^N g\left({x_{o(n)}}^\dagger \Sigma^{-1}x_{o(n)}\right)
\end{equation}

with $o$ being a permutation of $\left[\left[1;N\right]\right]$ such that:

\begin{displaymath}
g\left({x_{o(1)}}^\dagger \Sigma^{-1}x_{o(1)}\right) \leq \dots \leq g\left({x_{o(n)}}^\dagger \Sigma^{-1}x_{o(n)}\right)
\end{displaymath}

This therefore suggests to mix the partial estimation procedure with an M-estimation procedure, whenever $g'$ verifies the necessary conditions for convergence of the corresponding M-estimator \cite{marona1998robust}\cite{ollila2003robust} in a manner similar as the treatment given for Tyler's estimator.\\

Since these conditions include that $g'\geq 0$  for standard $M$-estimators of the covariance matrix \cite{marona1998robust}\cite{ollila2003robust}, this implies that $g$ is increasing and therefore the likelihood ordering of samples can again be done by ordering $\left({x_n}^\dagger R^{-1} x_n\right)_{1 \leq n \leq N}$. This thus corresponds to the following partial M-estimation procedure:\\

\begin{onehalfspace}
\begin{samepage}
\begin{algorithmic}[1]
	\Function{pM\_est}{$g',\left(x_n\right)_{1 \leq n \leq N}, p, \epsilon, K_\text{max}$}
		\State $\Sigma_{-1} \gets I$
		\For{$k$ \textbf{from} $1$ \textbf{to} $K_\text{max}$}
			\For{$n$ \textbf{from} $1$ \textbf{to} $N$}
				\State $\tau_n \gets {x_n}^\dagger \Sigma_{-1} x_n$
			\EndFor\vspace{1.5 mm}
			\State $o \gets \text{argsort}_{\uparrow}\left(\left(\tau_n\right)_{1 \leq n \leq N}\right)$\vspace{1.5 mm}
			\State $\Sigma \gets \frac{1}{\lceil pN \rceil} \sum\limits_{n=1}^{\lceil pN \rceil} g\left(\tau_{o_{(n)}}\right)x_{o(n)} {x_{o(n)}}^\dagger$\vspace{1.5 mm}
			\If{$\text{tr}\left(\left(\Sigma_{-1} \Sigma - I\right)^2\right) \leq \epsilon$}\vspace{1.5 mm}
				\State \textbf{break}	
			\Else
			\textbf{else do:}
				\State $\Sigma_{-1} \gets \Sigma^{-1}$
			\EndIf
		\EndFor
		\State \Return $\Sigma_{-1}$
	\EndFunction
\end{algorithmic}
\end{samepage}
\end{onehalfspace}
\vspace{2 mm}

Note that the partial SCM estimator is a special case of this estimator for which $g'(t)=1$.

This can also be extended using other estimators of the correlation matrix in order to solve correlation constrained problems. This corresponds to the following partial $M$-estimator procedure:\\

\begin{onehalfspace}
\begin{samepage}
\begin{algorithmic}[1]
	\Function{pM\_of}{$g',e,{x_n}_{1 \leq n \leq N},p,\epsilon,K_\text{max}$}
		\State $\Sigma_{-1} \gets I$
		\For{$k$ \textbf{from} $1$ \textbf{to} $K_\text{max}$}
			\For{$n$ \textbf{from} $1$ \textbf{to} $N$}
				\State $\tau_n \gets x_n \Sigma^{-1} x_n$
			\EndFor\vspace{1.5 mm}
			\State $o \gets \text{argsort}_{\uparrow}\left(\left(\tau_n\right)_{1 \leq n \leq N}\right)$ \vspace{1.5 mm}
			\State $\Sigma \gets e\left( \left(\sqrt{g'\left(\tau_{o(n)}\right)}x_{o(n)}\right)_{1 \leq n \leq \lceil pN \rceil}\right)$\vspace{1.5 mm}
			\If{ $\text{tr}\left({\left(\Sigma_{-1}\Sigma-I\right)}^2\right) \leq \epsilon$}\vspace{1.5 mm}
				\State \textbf{break}
			\Else
				\State $\Sigma_{-1} \gets \Sigma^{-1}$
			\EndIf
		\EndFor
		\State \Return $\Sigma_{-1}$
	\EndFunction
\end{algorithmic}
\end{samepage}
\end{onehalfspace}
\vspace{2 mm}

Whenever $g'$ is not positive, one can resort to the geodesic method given in \cite{CVC_BT}. This however cannot be extended to other known estimators for constrained problems:\\

\begin{onehalfspace}
\begin{samepage}
\begin{algorithmic}[1]
	\Function{pM\_exp\_cov}{$g,g',\left(x_n\right)_{1 \leq n \leq N},p,\epsilon,K_\text{max}$}
		\State $\Sigma_{\frac{1}{2}} \gets I$
		\State $\Sigma_{-\frac{1}{2}} \gets I$
		\For{$k$ \textbf{from} $1$ \textbf{to} $K_\text{max}$}
			\For{$n$ \textbf{from} $1$ \textbf{to} $N$}
				\State $\tau_n \gets x_n \Sigma^{-1} x_n$
			\EndFor\vspace{1.5 mm}
			\State $o \gets \text{argsort}_{\uparrow}\left(\left(g\left(\tau_n\right)\right)_{1 \leq n \leq N}\right)$\vspace{1.5 mm}
			\State $S \gets \frac{1}{\lceil pN \rceil}\sum\limits_{n=1}^{\lceil pN \rceil} g'\left(\tau_{o(n)}\right)x_{o(n)}{x_{o(n)}}^\dagger$\vspace{1.5 mm}
			\State $\Sigma \gets \Sigma_{\frac{1}{2}}\exp\left(\Sigma_{-\frac{1}{2}} S \Sigma_{-\frac{1}{2}}-I\right)\Sigma_{\frac{1}{2}}$\vspace{1.5 mm}
			\If{$\text{tr}\left(\left(\Sigma_{-\frac{1}{2}} \Sigma \Sigma_{-\frac{1}{2}} - I\right)^2\right) \leq \epsilon$}\vspace{1.5 mm} 
				\State \textbf{break}	
			\Else	
				\State $\Sigma_{\frac{1}{2}} \gets \sqrt{\Sigma}$
				\State $\Sigma_{-\frac{1}{2}} \gets {\Sigma_{\frac{1}{2}}}^{-1}$
			\EndIf
		\EndFor
		\State \Return $\Sigma_{-\frac{1}{2}}$
	\EndFunction
\end{algorithmic}
\end{samepage}
\end{onehalfspace}
\vspace{2 mm}

For example for the case of a gaussian distribution in $\mathbb{C}^d$ and under cirularity hypothesis, one has \cite{CVC_BT}:

\begin{displaymath}
g(t) = t-\frac{1}{2}\log(t)
\end{displaymath}

Thus one can adapt the geodesic partial M estimation procedure as follows: 

\begin{onehalfspace}
\begin{samepage}
\begin{algorithmic}[1]
	\Function{pcg\_cov}{$\left(x_n\right)_{1 \leq n \leq N},\epsilon,K_\text{max}$}\vspace{1.5 mm}
		\State $\Sigma \gets \frac{1}{N}\sum\limits_{n=1}^{N} x_n {x_n}^\dagger$\vspace{1.5 mm}
		\State $\Sigma_{\frac{1}{2}} \gets \sqrt{\Sigma}$
		\State $\Sigma_{-\frac{1}{2}} \gets {\Sigma_{\frac{1}{2}}}^{-1}$
		\For{$k$ \textbf{from} $1$ \textbf{to} $K_\text{max}$}
			\For{$n$ \textbf{from} $1$ \textbf{to} $N$}
				\State $\tau_n \gets {x_n}^\dagger \Sigma^{-1} x_n$
			\EndFor\vspace{1.5 mm}
			\State $o \gets \text{argsort}_{\uparrow}\left(\left(\tau_n-\frac{1}{2}\log\left(\tau_n\right)\right)_{1 \leq n \leq N}\right)$\vspace{1.5 mm}
			\State $S \gets \frac{1}{\left(1-\frac{1}{2d}\right)\lceil pN \rceil}\sum\limits_{n=1}^{\lceil pN \rceil} \left(1-\frac{1}{2\tau_{o(n)}}\right)x_{o(n)}{x_{o(n)}}^\dagger$\vspace{1.5 mm}
			\State $R \gets \Sigma_{\frac{1}{2}}\exp\left(\Sigma_{-\frac{1}{2}} S \Sigma_{-\frac{1}{2}}-I\right)\Sigma_{\frac{1}{2}}$\vspace{1.5 mm}
			\If{$\text{tr}\left(\left(\Sigma_{-\frac{1}{2}} \Sigma \Sigma_{-\frac{1}{2}} - I\right)^2\right) \leq \epsilon$}\vspace{1.5 mm}
				\State \textbf{break}	
			\Else	
				\State $\Sigma_{\frac{1}{2}} \gets \sqrt{\Sigma}$
				\State $\Sigma_{-\frac{1}{2}} \gets {\Sigma_{\frac{1}{2}}}^{-1}$
			\EndIf
		\EndFor
		\State \Return $\Sigma_{-\frac{1}{2}}$
	\EndFunction
\end{algorithmic}
\end{samepage}
\end{onehalfspace}
\vspace{2 mm}

\section{Simulations}

We now show some simulation results of adaptive detectors using various estimators introduced in this article, using the detection tests introduced in \cite{CVC_BT}\cite{scharf1994matched}.\\

The simulation results are shown in a single channel scenario of dimension $d=8$, as a function of the SiNR of the target signal.\\
the background noise is generated as a white gaussian noise of unit variance.\\
The target signal is generated as a complex centered circular 1-dimensional gaussian signal aligned with the test signal $s$, whose variance $\sigma$ is such that:

\begin{displaymath}
10 \log_{10}(\sigma) = \text{SiNR}
\end{displaymath}

The detection thresholds are defined to have a false alarm rate of $10^{-4}$; they are learned on clean training sets, that is that they contain no outliers.\\

\subsection{Impact of outliers on adaptive detection}

In this scenario, outlier target signals are present among the $N$ samples used in the estimation of the prior covariance, following the same law as the target signal but being independently drawn.\\

Prior covariances are estimated with $N = 22$ independently drawn noise samples, among which a varying number of outliers is present, going from 0 to 5 in order to show their impact on the detection performances of adaptive filters.\\

\begin{figure}[H]
	\includegraphics[width=\linewidth]{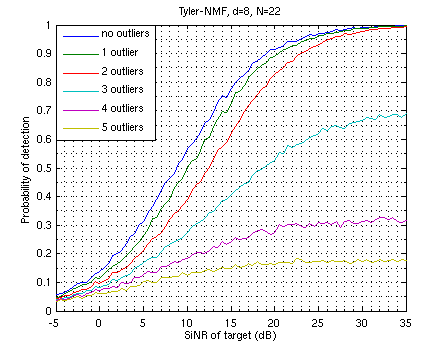}
	\caption{Detection capability Tyler-NMF under the influence of secondary target signals}
	\label{fig_Tyler_outlier}
\end{figure}

\begin{figure}[H]
	\includegraphics[width=\linewidth]{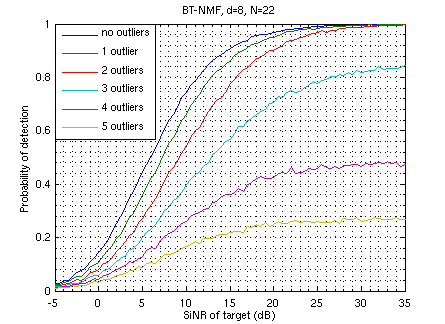}
	\caption{Detection capability BT-NMF under the influence of secondary target signals}
	\label{fig_BT_outlier}
\end{figure}

\begin{figure}[H]
	\includegraphics[width=\linewidth]{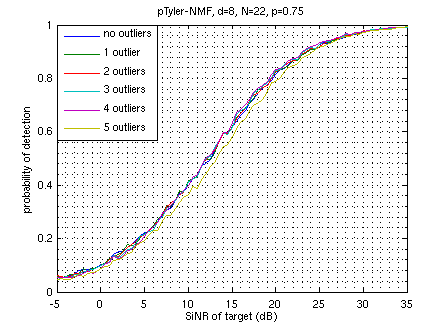}
	\caption{Detection capability pTyler-NMF under the influence of secondary target signals}
	\label{fig_pTyler_outlier}
\end{figure}

\begin{figure}[H]
	\includegraphics[width=\linewidth]{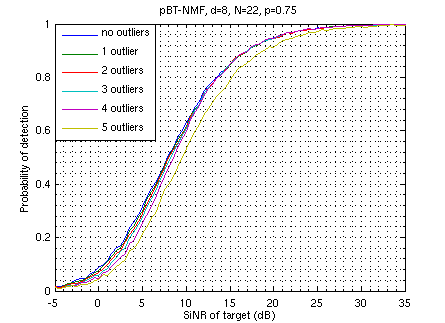}
	\caption{Detection capability pBT-NMF under the influence of secondary target signals}
	\label{fig_pBT_outlier}
\end{figure}

Figures~\ref{fig_Tyler_outlier},~\ref{fig_BT_outlier},~\ref{fig_pTyler_outlier} and ~\ref{fig_pBT_outlier}  respectively show the detection performances of the NMF test using the Tyler, BT, pTyler and pBT estimators, with a partial order of p = 0.75. As is visible here, the partial estimation almost completely mitigates the outliers' impact on th detection performances for Tyler-type estimators.\\

\begin{figure}[H]
	\includegraphics[width=\linewidth]{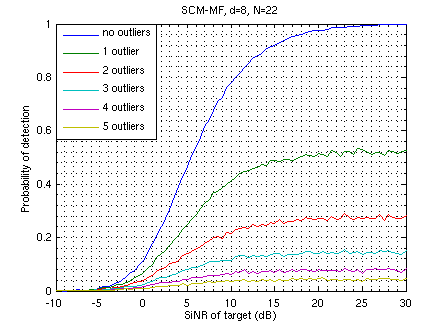}
	\caption{Detection capability SCM-MF under the influence of secondary target signals}
	\label{fig_SCM_outlier}
\end{figure}

\begin{figure}[H]
	\includegraphics[width=\linewidth]{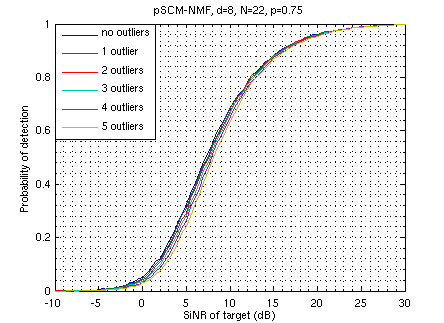}
	\caption{Detection capability pSCM-MF under the influence of secondary target signals}
	\label{fig_pSCM_outlier}
\end{figure}

Figures~\ref{fig_SCM_outlier} and~\ref{fig_pSCM_outlier} show the detection performances of the matched filter test (MF) using the SCM and pSCM estimators, with a partial order of p = 0.75. Again the partial estimation almost completely mitigates the outliers' impact on the detection performances in this case.

\begin{figure}[H]
	\includegraphics[width=\linewidth]{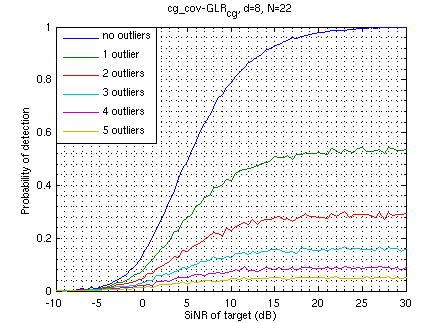}
	\caption{Detection capability cg\_cov-$\text{GLR}_\text{cg}$ under the influence of secondary target signals}
	\label{fig_cg_outlier}
\end{figure}

\begin{figure}[H]
	\includegraphics[width=\linewidth]{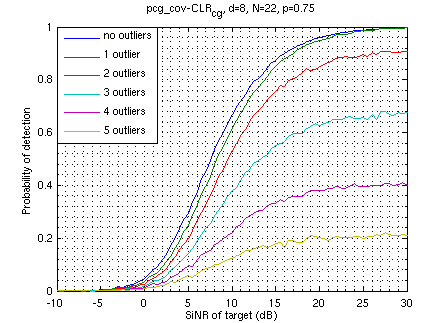}
	\caption{Detection capability pcg\_cov-$\text{GLR}_\text{cg}$ under the influence of secondary target signals}
	\label{fig_pcg_outlier}
\end{figure}

Finally figures~\ref{fig_cg_outlier} and~\ref{fig_pcg_outlier} show the detection performances of the $\text{GLR}\_\text{cg}$ test using the cg\_cov and pcg\_cov estimators, with a partial order of p = 0.75. Again the effect of the outliers is mitigated, although their impact is still seen when more than 3 of them are present among the N= 22 samples used for the estimation.

Finally, one should note that although partial estimators allow to mitigate the impact of outliers, they also degrade the optimal detection performances. This is understandable as fewer samples are actually used to compute the estimates. Thus one has to carefully balance the order of the partial estimation with the desired performances, depending on applications.

\section{Conclusion}

We have presented a theoretical background for a likelihood based data selection scheme, in order to solve problems related to outlier detection and rejection in estimation procedures. This leads to the notion of partial estimators, which can be used in order to produce new estimation procedures which are robust to outliers from known estimators.\\
This principle is then applied to several covariance estimators, such as the usual sample covariance matrix, Tyler's estimator and other type of M-estimators, as well as the Burg-Tyler estimator \cite{CVC_BT}.\\
The authors would like to outline the particularly important fact that such estimation procedures are best used in cases in which there exists a reference measure of the underlying space; indeed otherwise such a reference measure has to be specified, which then creates a bias towards the chosen reference. 

\ifCLASSOPTIONcaptionsoff
  \newpage
\fi



%

\bibliographystyle{plain}
\bibliography{CES_biblio}

\begin{thebibliography}{10}

\bibitem{akaike1998information}
Hirotogu Akaike.
\newblock Information theory and an extension of the maximum likelihood
  principle.
\newblock In {\em Selected Papers of Hirotugu Akaike}, pages 199--213.
  Springer, 1998.

\bibitem{anderson1985maximum}
Theodore~Wilbur Anderson and Ingram Olkin.
\newblock Maximum-likelihood estimation of the parameters of a multivariate
  normal distribution.
\newblock {\em Linear algebra and its applications}, 70:147--171, 1985.

\bibitem{aubry2013median}
A~Aubry, A~De~Maio, Luca Pallotta, A~Farina, and C~Fantacci.
\newblock Median matrices and geometric barycenters for training data
  selection.
\newblock In {\em Radar Symposium (IRS), 2013 14th International}, volume~1,
  pages 331--336. IEEE, 2013.

\bibitem{barbaresco1997analyse}
Fr{\'e}d{\'e}ric BARBARESCO.
\newblock Analyse spectrale par decomposition recursive en sous-espaces propres
  via les coefficients de reflexion.
\newblock In {\em 16° Colloque sur le traitement du signal et des images, FRA,
  1997}. GRETSI, Groupe d’Etudes du Traitement du Signal et des Images, 1997.

\bibitem{burg1978maximum}
JP~Burg.
\newblock Maximum entropy spectral analysis, modern spectrum analysis dg
  childers, 34--41, 1978.

\bibitem{chmielewski1981elliptically}
MA~Chmielewski.
\newblock Elliptically symmetric distributions: a review and bibliography.
\newblock {\em International Statistical Review/Revue Internationale de
  Statistique}, pages 67--74, 1981.

\bibitem{conte1991modelling}
E~Conte, M~Longo, and M~Lops.
\newblock Modelling and simulation of non-rayleigh radar clutter.
\newblock In {\em Radar and Signal Processing, IEE Proceedings F}, volume 138,
  pages 121--130. IET, 1991.

\bibitem{conte1987characterisation}
Ernesto Conte and Maurizio Longo.
\newblock Characterisation of radar clutter as a spherically invariant random
  process.
\newblock {\em Communications, Radar and Signal Processing, IEE Proceedings F},
  134(2):191--197, 1987.

\bibitem{conte1995asymptotically}
Ernesto Conte, Marco Lops, and Giuseppe Ricci.
\newblock Asymptotically optimum radar detection in compound-gaussian clutter.
\newblock {\em Aerospace and Electronic Systems, IEEE Transactions on},
  31(2):617--625, 1995.

\bibitem{cousineau2015outliers}
Denis Cousineau and Sylvain Chartier.
\newblock Outliers detection and treatment: a review.
\newblock {\em International Journal of Psychological Research}, 3(1):58--67,
  2015.

\bibitem{CVC_BT}
Christophe Culan and Claude Adnet.
\newblock {Maximum likelihood estimation of covariances of elliptically
  symmetric distributions}.
\newblock {\em eprint}, arXiv:1611.04365, 2016.

\bibitem{daszykowski2007robust}
Michal Daszykowski, Krzysztof Kaczmarek, Yvan Vander~Heyden, and Beata Walczak.
\newblock Robust statistics in data analysis—a review: basic concepts.
\newblock {\em Chemometrics and intelligent laboratory systems},
  85(2):203--219, 2007.

\bibitem{decurninge2014burg}
Alexis Decurninge and Frederic Barbaresco.
\newblock Burg estimation of radar covariance matrix for mixtures of gaussian
  stationary distributions.
\newblock In {\em Radar Conference (Radar), 2014 International}, pages 1--6.
  IEEE, 2014.

\bibitem{decurninge2016robust}
Alexis Decurninge and Fr{\'e}d{\'e}ric Barbaresco.
\newblock Robust burg estimation of radar scatter matrix for mixtures of
  gaussian stationary autoregressive vectors.
\newblock {\em arXiv preprint arXiv:1601.02804}, 2016.

\bibitem{gini2002vector}
Fulvio Gini and Alfonso Farina.
\newblock Vector subspace detection in compound-gaussian clutter. part i:
  survey and new results.
\newblock {\em Aerospace and Electronic Systems, IEEE Transactions on},
  38(4):1295--1311, 2002.

\bibitem{haykin1982maximum}
Simon Haykin, Brian~W Currie, and Stanislav~B Kesler.
\newblock Maximum-entropy spectral analysis of radar clutter.
\newblock {\em Proceedings of the IEEE}, 70(9):953--962, 1982.

\bibitem{krishnaiah1986complex}
PR~Krishnaiah and Jugan Lin.
\newblock Complex elliptically symmetric distributions.
\newblock {\em Communications in Statistics-Theory and Methods},
  15(12):3693--3718, 1986.

\bibitem{kullback1951information}
Solomon Kullback and Richard~A Leibler.
\newblock On information and sufficiency.
\newblock {\em The annals of mathematical statistics}, 22(1):79--86, 1951.

\bibitem{marona1998robust}
RA~Marona and VJ~Yohai.
\newblock Robust estimation of multivariate location and scatter.
\newblock {\em Kotz, S., Read, Banks, D.(Eds.), Encyclopedia of Statistical
  Sciences}, 2:590, 1998.

\bibitem{ollila2003robust}
Esa Ollila and Visa Koivunen.
\newblock Robust antenna array processing using m-estimators of
  pseudo-covariance.
\newblock In {\em Personal, Indoor and Mobile Radio Communications, 2003. PIMRC
  2003. 14th IEEE Proceedings on}, volume~3, pages 2659--2663. IEEE, 2003.

\bibitem{ollila2012complex}
Esa Ollila, David~E Tyler, Visa Koivunen, and H~Vincent Poor.
\newblock Complex elliptically symmetric distributions: Survey, new results and
  applications.
\newblock {\em Signal Processing, IEEE Transactions on}, 60(11):5597--5625,
  2012.

\bibitem{pascal2008covariance}
Fr{\'e}d{\'e}ric Pascal, Yacine Chitour, Jean-Philippe Ovarlez, Philippe
  Forster, and Pascal Larzabal.
\newblock Covariance structure maximum-likelihood estimates in compound
  gaussian noise: Existence and algorithm analysis.
\newblock {\em Signal Processing, IEEE Transactions on}, 56(1):34--48, 2008.

\bibitem{pena2012multivariate}
Daniel Pe{\~n}a and Francisco~J Prieto.
\newblock Multivariate outlier detection and robust covariance matrix
  estimation.
\newblock {\em Technometrics}, 2012.

\bibitem{scharf1991statistical}
Louis~L Scharf.
\newblock {\em Statistical signal processing}, volume~98.
\newblock Addison-Wesley Reading, MA, 1991.

\bibitem{scharf1994matched}
Louis~L Scharf and Benjamin Friedlander.
\newblock Matched subspace detectors.
\newblock {\em Signal Processing, IEEE Transactions on}, 42(8):2146--2157,
  1994.

\bibitem{trench1964algorithm}
William~F Trench.
\newblock An algorithm for the inversion of finite toeplitz matrices.
\newblock {\em Journal of the Society for Industrial and Applied Mathematics},
  12(3):515--522, 1964.

\bibitem{tyler1987statistical}
David~E Tyler.
\newblock Statistical analysis for the angular central gaussian distribution on
  the sphere.
\newblock {\em Biometrika}, 74(3):579--589, 1987.

\bibitem{ulrych1976time}
Tad~J Ulrych and Rob~W Clayton.
\newblock Time series modelling and maximum entropy.
\newblock {\em Physics of the Earth and Planetary Interiors}, 12(2):188--200,
  1976.

\bibitem{zohar1969toeplitz}
Shalhav Zohar.
\newblock Toeplitz matrix inversion: The algorithm of wf trench.
\newblock {\em Journal of the ACM (JACM)}, 16(4):592--601, 1969.

\end{thebibliography}



%

\begin{IEEEbiographynophoto}{Christophe Culan}
	was born in Villeneuve-St-Georges, France, on November $\text{1}^\text{st}$, 1988. He received jointly the engineering degree from Ecole Centrale de Paris (ECP), France and a master's degree in Physico-informatics from Keio University, Japan, in 2013.\\
	He was a researcher in applied physics in Itoh laboratory from 2013 to 2014, specialized in quantum information and quantum computing, and has contributed to several publications related to these subjects.\\
	He currently holds a position as a research engineer in Thales Air Systems, Limours, France, in the Advanced Radar Concepts division. His current research interests include statistical signal and data processing, robust statistics, machine learning and information geometry.
\end{IEEEbiographynophoto}

\begin{IEEEbiographynophoto}{Claude Adnet}
	was born in Aÿ, France, in 1961. He received the DEA degree in signal processing and Phd degree in 1988 and 1991 respectively, from the Institut National Polytechnique de Grenoble (INPG), Grenoble France.
	Since then, he has been working for THALES Group, where he is now Senior Scientist. His research interests include radar signal  processing and radar data processing.
\end{IEEEbiographynophoto}






\end{document}